\documentclass[preprint,12pt]{elsarticle}

\usepackage{amssymb}
\usepackage{amscd,amsmath,amsfonts,graphicx}
\usepackage{amsthm}

\begin{document}

\baselineskip=22pt \centerline{\Large \bf The Bivariate
Lack-of-Memory Distributions}  \vspace{1cm} \centerline{Gwo Dong
Lin, Xiaoling Dou and Satoshi Kuriki} \centerline{Academia Sinica,
Taiwan, Waseda University, Japan}  \centerline{and The Institute of
Statistical Mathematics, Japan} \vspace{1cm} \noindent {\bf
Abstract.} We treat all the bivariate lack-of-memory (BLM)
distributions in a unified approach and develop some new general
properties of the BLM distributions, including joint moment
generating function, product moments and dependence structure.
Necessary and sufficient conditions for the survival functions of
BLM distributions to be totally positive of order two are given.
Some previous results about specific BLM distributions are improved.
In particular, we show that both the Marshall--Olkin survival copula
and survival function are totally positive of all orders, regardless
of parameters. Besides, we point out that Slepian's inequality  also
holds true for BLM distributions.

\vspace{0.7cm} \hrule \medskip \vspace{0.2cm}\noindent{\bf 2010 AMS
Mathematics Subject Classifications}. Primary: 62N05, 62N86, 62H20, 62H05.
{\bf Postal addresses:} Gwo Dong Lin, Institute of Statistical
Science, Academia Sinica, Taipei
11529, Taiwan, R.O.C. (E-mail: gdlin@stat.sinica.edu.tw)\\
Xiaoling Dou, Waseda University,
3-4-1 Ohkubo, Shinjuku, Tokyo 169-8555, Japan (E-mail: xiaoling@aoni.waseda.jp)\\
Satoshi Kuriki, The Institute of Statistical Mathematics, 10-3
Midoricho, Tachikawa, Tokyo 190-8562, Japan  (E-mail:
kuriki@ism.ac.jp)

\newpage \noindent {\bf 1. Introduction}

The classical univariate lack-of-memory (LM) property is a
remarkable characterization of the exponential distribution which
plays a prominent role in reliability theory, queuing theory and
other applied fields (Feller 1965, Fortet 1977, Galambos and Kotz
1978). The recent bivariate LM property is, however, shared by
 the famous Marshall and Olkin's, Block
and Basu's as well as Freund's bivariate exponential distributions,
among many others; see, e.g., Chapter 10 of Balakrishnan and Lai
(2009), Chapter 47 of Kotz et al.\,(2000) and Kulkarni (2006). These
bivariate distributions have been well investigated individually in
the literature. Our main purpose in this paper is, however, to
develop in a {unified approach} some new general properties of  the
bivariate lack-of-memory (BLM) distributions which share the same
bivariate LM property.

In Section 2, we first review the univariate and bivariate LM
properties, and then summarize the important known properties of the
BLM distributions. We derive in Section 3 some new general
properties of the BLM distributions, including joint moment
generating function, product moments and stochastic inequalities.
The dependence structures of the BLM distributions are investigated
in Section 4. We find necessary and sufficient conditions for the
survival functions (and the densities if they exist) of BLM
distributions to be totally positive of order two. Some previous
results about specific BLM distributions are improved. In
particular, we show that both the Marshall--Olkin survival copula
and survival function are totally positive of all orders, regardless
of parameters. In Section 5, we study the stochastic comparisons in
the family of all
 BLM distributions and point out that Slepian's lemma/inequality for
bivariate normal distributions also holds true for BLM
distributions.
\medskip\\
\noindent {\bf 2.  Lack-of-Memory Property}

We first review the well-known univariate lack-of-memory property.
Let $X$ be a nonnegative random variable with distribution function
$F.$ Then $F$ satisfies (multiplicative) Cauchy's functional
equation
\begin{eqnarray}
{\overline{F}(x+y)=\overline{F}(x)\overline{F}(y), \ x\ge 0, y\ge
0,}
\end{eqnarray} where ${\overline{F}(x)=1-F(x)=\Pr(X>x)},$
if and only if $F(0)=1$ ($X$ degenerates at 0) or
$F(x)=1-\exp({-\lambda x}),\ x\ge 0,$ for some constant $\lambda
>0,$ denoted by $X\sim Exp(\lambda)$ ($X$ has an exponential
distribution with positive parameter $\lambda$).
 If $X$ is the lifetime of a system with positive survival
 function
$\overline{F},$ then Eq.\,(1) is equivalent to \begin{eqnarray}{\Pr
(X>x+y|\ X>y)=\Pr(X>x),~~x\ge 0, \ y\ge 0.}\end{eqnarray} This means
that the conditional probability of a system surviving to time $x+y$
given surviving to time $y$ is equal to the unconditional
probability of the system surviving to time $x$. Namely, the failure
performance of the system does not depend on the past, given its
present condition. In such a case (2), we say that the distribution
$F$ lacks memory at each point $y.$ So Eq.\,(1) is called the LM
property or memoryless property of $F.$

For simplicity, we consider only {\it positive} random variable
$X\sim F$  from now on. Then, the LM property (1) holds true iff
$X\sim Exp(\lambda)$ for some $\lambda>0.$

We next consider the bivariate LM property. Let the positive random
variables $X$ and $Y$ have joint distribution $H$ with marginals $F$
and $G$. Namely, $(X,Y)\sim H,\ X\sim F,\ Y\sim G$. Moreover, denote
the survival function of $H$ by
\begin{eqnarray*}\overline{H}(x,y)\equiv \Pr(X>x,\,Y>y)
=1-F(x)-G(y)+H(x,y),\  x,y\ge 0.
\end{eqnarray*}
An intuitive extension of the LM property (2) to the bivariate case
is the strict BLM property:
$${\Pr(X>x+s,\,Y>y+t|\ X>s,Y>t)=\Pr(X>x,\,Y>y),\ x,y,s,t\ge 0}$$
({$H$ lacks memory at each pair $(s,t)$}), which is equivalent to
\begin{eqnarray}\overline{H}(x+s,\,y+t)=\overline{H}(x,y)\overline{H}(s,t),\
\forall\ x,y,s,t\ge 0,
\end{eqnarray} if the survival function $\overline{H}$ is
positive.  In a two-component system, this means as before that the
conditional probability of two components surviving to times
$(x+s,y+t)$ given surviving to times $(s, t)$ is equal to the
unconditional probability of these two components surviving to times
$(x,y)$. But Eq.\,(3) has only one solution (Marshall and Olkin
1967, p.\,33), namely, the independent bivariate exponential
distribution with survival function
$${\overline{H}(x,y)=\exp[{-(\lambda x+\delta y)}], \ x,y\ge
0,}$$ for some constants $\lambda,\delta>0;$ in other words, $X$ and
$Y$ are independent random variables and  $X\sim Exp(\lambda),\
Y\sim  Exp(\delta)$ for some positive parameters $\lambda, \delta$.

 In their pioneering paper, Marshall and Olkin (1967) considered instead the weaker BLM property
({with $s=t$})
$$\Pr(X>x+t,\,Y>y+t|\ X>t,\,Y>t)=\Pr(X>x,\,Y>y),\ x,y,t\ge 0$$
({$H$ lacks memory at each {\it equal} pair $(t,t)$}), and solved
the functional equation
\begin{eqnarray}{\overline{H}(x+t,\,y+t)=\overline{H}(x,y)\overline{H}(t,t),\
\forall\ x,y,t\ge0.} \end{eqnarray} It turns out that for given
$(X,Y)\sim H$ with marginals $F, G$ on $(0,\infty)$, $H$ satisfies
the BLM property (4) iff its survival function is of the form
\begin{equation}
   { \overline{H}(x,y)=\left\{\begin{array}{cc}
                    e^{-\theta y}\,\overline{F}(x-y), & x\ge y\ge 0 \vspace{0.1cm}\\
                    e^{-\theta x}\,\overline{G}(y-x),  & y\ge x\ge 0,
                  \end{array}\right.}
\end{equation}
where $\theta$  is a positive constant (see also Barlow and Proschan
1981, p.\,130).

 For convenience, denote by $BLM(F,G,\theta)$ the BLM distribution
 $H$ with marginals $F,G$, parameter $\theta>0$ and
survival function $\overline{H}$ in (5), and denote by ${\cal BLM}$
the family of all BLM distributions, namely,
$${\cal BLM}=\{H: H=BLM(F,G,\theta),\ \hbox{where}\ \theta>0,\
\hbox{and}\ F,\ G\ \hbox{are marginal distributions} \}.$$  Theorem
1 below summarizes some important known properties of the BLM
distributions; for more details, see Marshall and Olkin (1967),
Block and Basu (1974), Block (1977), and Ghurye and Marshall (1984).
For convenience, denote $a\vee b=\max\{a,b\}$ and $a\wedge
b=\min\{a,b\}.$\\ \noindent{\bf Theorem 1.} Let $(X,Y)\sim
H=BLM(F,G,\theta)\in{\cal BLM}.$
Then the following statements are true.\\
 (i) The
marginals $F$, $G$ have densities $f$, $g$, respectively. Moreover,
the right-hand derivatives\\ $f(x)=\lim_{\varepsilon\to
0^+}[F(x+\varepsilon)-F(x)]/{\varepsilon}$\ and
$g(x)=\lim_{\varepsilon\to
0^+}[G(x+\varepsilon)-G(x)]/{\varepsilon}$ exist\ for all $x\ge 0$,
which are right-continuous and are of bounded variation  on $[0,\infty)$.\\
(ii) $\Pr(X-Y>t)=\overline{F}(t)-f(t)/\theta$\ and
$\Pr(Y-X>t)=\overline{G}(t)-g(t)/\theta$ for all $t\ge 0.$\\
(iii) Both $e^{\theta x}f(x)$ and $e^{\theta x}g(x)$ are
increasing (nondecreasing) in $x\ge 0.$ \\
(iv) $F(x)+G(x)\ge 1-\exp({-\theta x}),\ x\ge 0.$\\
(v)  $X\wedge Y\sim
Exp(\theta)$ and is independent of $X-Y.$ \\
(vi) $f(0)\vee g(0)\le \theta\le f(0)+g(0).$ \\
(vii) $f^{\prime}(x)+\theta f(x)\ge 0,\ g^{\prime}(x)+\theta g(x)\ge
0,\ x\ge 0,$\ \ if $f$ and $g$ are differentiable.
\smallskip\\
\noindent{\bf Remark 1.} Some of the above necessary conditions
(i)--(vii) also play as sufficient conditions for $(X,Y)$ to obey a
BLM distribution. For example, in addition to the above conditions
(vi) and (vii), assume that the marginal densities are  absolutely
continuous, then the $\overline{H}$ in (5) is a {\it bona fide}
survival function. This is a slight modification of Theorem 5.1 of
Marshall and Olkin (1967) who required (vi$^{\prime}$)
$[f(0)+g(0)]/2\le \theta\le f(0)+g(0)$ instead of (vi) above. Note
that conditions (vi) and (vi$^{\prime}$) are different unless
$f(0)=g(0)$, and that (vi) is a consequence of (iii) and (iv) (see
Corollary 2(i) below and Ghurye and Marshall 1984, p.\,792). On the
other hand, the condition (v) together with continuous marginals
$F,G$ also implies that $(X,Y)$ has a BLM distribution (Block 1977,
p.\,810). It is interesting to recall that for {\it independent}
nondegenerate random variables $X$ and $Y$, the above independence
of $X\wedge Y$ and $X-Y$ is a characterization of the
exponential/geometric distributions under suitable conditions (see
Ferguson  1964, 1965, Crawford 1966, and Rao and Shanbhag 1994).
Namely, in general, the BLM distributions share the same
independence property of $X\wedge Y$ and $X-Y$ with independent
exponential/geometric random variables.
\smallskip\\
\noindent{\bf Remark 2.}  There are some more observations: (a)
$\Pr(X=Y)=[f(0)+g(0)]/\theta-1$ by the above (ii), (b) at least one
of $f(0)$ and $g(0)$ is positive, (c) the survival function
$\overline{H}$ in (5) is purely singular (i.e., $X=Y$ almost surely)
iff $\theta=[f(0)+g(0)]/2$ iff $f(0)=g(0)=\theta$ (because $f(0)\ne
g(0)$ implies $\theta>[f(0)+g(0)]/2$\ by (vi)), and (d)
$\overline{H}$
 is absolutely continuous (i.e.,
$X\ne Y$ almost surely) iff the marginal densities together satisfy
$f(0)+g(0)=\theta$ (see Ghurye and Marshall 1984, p.\,792). In view
of the above results, the survival function (5) of
$H=BLM(F,G,\theta)$ can be rewritten as the convex combination of
two extreme ones:
$${\overline{H}(x,y)}=
\left(2-\frac{f(0)+g(0)}{\theta}\right){\overline{H}{\!}_a(x,y)}+
\left(\frac{f(0)+g(0)}{\theta}-1\right){\overline{H}{\!}_s(x,y)} ,\
x,y\ge 0,$$ where  $H_a$ is absolutely continuous and $H_s$ is
purely singular with survival function
$\overline{H}{\!}_s(x,y)=\exp[-\theta\max\{x,y\}],\,x,y\ge 0.$
Clearly, the parameter $\theta$ regulates $H$ between $H_a$ and
$H_s.$  On the other hand,  Ghurye and Marshall (1984, Section 3)
gave an interesting random decomposition of $(X,Y)\sim H\in{\cal
BLM}$ and represented $\overline{H}$ as a Laplace--Stieltjes
integral by another bivariate survival function. See also Ghurye
(1987) and Marshall and Olkin (2015) for further generalizations of
the BLM distributions.
\smallskip\\
\noindent{\bf Remark 3.} Kulkarni (2006) proposed an interesting and
useful approach to construct some BLM distributions by starting with
marginal failure rate functions.
 First, choose two real-valued functions $r_1, r_2$ and a constant $\theta$ satisfying the following (modified) conditions:\\
 (a) The functions $r_i, i=1,2$, are {\it absolutely continuous}  on $[0,\infty)$ and $\theta>0$.\\
 (b) $0\le r_i(x)\le \theta,\ x\ge 0,\ i=1,2.$\\
 (c) $\int_0^{\infty}r_i(x)dx=\infty,\ i=1,2.$\\
 (d) $r_i(x)(\theta-r_i(x))+r_i^{\prime}(x)\ge 0,\ x\ge 0,\ i=1,2.$\\
 (e) $r_1(0)+r_2(0)\ge \theta.$\\
 Then set {$\overline{F}(x)=\exp(-\int_0^xr_1(t)dt), x\ge 0,$ and $\overline{G}(y)=\exp(-\int_0^yr_2(t)dt), y\ge
 0.$} In this way, the $H$ defined through (5) is a {\it bona fide} BLM
 distribution because the above conditions (a)--(e) together imply  that conditions (vi) and (vii) in Theorem 1 hold true (see Remark 1).
Conversely, under the above smoothness conditions on $r_i$ and the
setting of $F$ and $G$, if $H$ in (5) is a BLM distribution, then
its marginal failure rate functions $r_i$ should satisfy conditions
(b)--(e) from which some properties in Theorem 1 follow immediately
(Kulkarni 2006, Proposition 1).

We now recall three important BLM distributions in the literature.
For more details, see, e.g.,  Chapter 10 of Balakrishnan and Lai
(2009).
\smallskip\\
\noindent{\bf Example 1.} Marshall and Olkin's (1967) bivariate
exponential distribution (BVE)\\ If {both marginals $F$ and $G$ are
exponential}, then the $BLM(F,G,\theta)\in{\cal BLM}$ defined in (5)
reduces to the Marshall--Olkin BVE with survival function of the
form
\begin{eqnarray}{\overline{H}(x,y)}&=&{\exp[-\lambda_1
x-\lambda_2 y-\lambda_{12}\max\{x,y\}]}\\
&\equiv& \frac{\lambda_1+\lambda_2}{\lambda}{\overline{H}{\!}_a(x,y)}+
\frac{\lambda_{12}}{\lambda}{\overline{H}{\!}_s(x,y)} ,\ \ x,y\ge 0,
\end{eqnarray} where $\lambda_1,\lambda_2, \lambda_{12}$ are
positive constants, $\lambda=\lambda_1+\lambda_2+\lambda_{12}$, and
$H_a$, $H_s$ (written explicitly below) are absolutely continuous
and singular bivariate distributions, respectively.

In practice, the Marshall--Olkin BVE arises from a shock model for
a two-component system. Formally,   the lifetimes of two components
are $(X, Y)=(X_1\wedge X_3,\ X_2\wedge X_3),$ where $X_1\sim
Exp(\lambda_1),$ $X_2\sim
   Exp(\lambda_2)$ and $X_3\sim  Exp(\lambda_{12})$ are independent.
So they have a joint survival function $\overline{H}$ defined in
(6). The singular part in (7) is identified by the conditional
probability: $\overline{H}{\!}_s(x,y)=\Pr(X>x,\,Y>y|\ X_3\le
X_1\wedge X_2)=\exp[-\lambda\max\{x,y\}],$ while the absolutely
continuous part $\overline{H}{\!}_a$ is calculated from
$\overline{H}$ and $\overline{H}{\!}_s$ via (7) (see the next
example).
\smallskip\\
\noindent{\bf Example 2.}  Block and Basu's (1974) bivariate
exponential
distribution\\
The Block--Basu BVE is actually the absolute continuous part $H_a$
of Marshall--Olkin BVE in (7) and has a joint density of the form
\begin{eqnarray}
   { {h}(x,y)=\left\{\begin{array}{cc}
                    \frac{\lambda_2\lambda(\lambda_1+\lambda_{12})}{\lambda_1+\lambda_2}
                    \exp[{-(\lambda_1+\lambda_{12})x-\lambda_2y}], &
                    x\ge
                    y> 0\vspace{0.2cm}\\
                    \frac{\lambda_1\lambda(\lambda_2+\lambda_{12})}{\lambda_1+\lambda_2}
                    \exp[{-\lambda_1x-(\lambda_2+\lambda_{12})y}], &
                    y>
                    x> 0,
                  \end{array}\right.}
\end{eqnarray}
where $\lambda_1,\lambda_2,\lambda_{12}>0,$ and
$\lambda=\lambda_1+\lambda_2+\lambda_{12}.$ Its survival function is
equal to
\begin{eqnarray*}& &\overline{H}(x,y)={\overline{H}{\!}_a(x,y)}\\
&=&\frac{\lambda}{\lambda_1+\lambda_2}\exp[-\lambda_1 x-\lambda_2
y-\lambda_{12}\max\{x,y\}]
-\frac{\lambda_{12}}{\lambda_1+\lambda_2}\exp[-\lambda\max\{x,y\}],\
x,y\ge 0.
\end{eqnarray*}

Note that in this case, the  marginals are not exponential but
rather  {\it negative mixtures} of two exponentials. Specifically,
$\overline{F}(x)=\frac{\lambda}{\lambda_1+\lambda_2}\exp[-(\lambda_1+\lambda_{12})x]-
\frac{\lambda_{12}}{\lambda_1+\lambda_2}\exp({-\lambda x}),\ x\ge
0,$ and
$\overline{G}(y)=\frac{\lambda}{\lambda_1+\lambda_2}\exp[{-(\lambda_2+\lambda_{12})y}]-
\frac{\lambda_{12}}{\lambda_1+\lambda_2}\exp({-\lambda y}),\ y\ge
0.$
\smallskip\\
\noindent{\bf Example 3.}  Freund's (1961) bivariate exponential
distribution\\
The Freund BVE has a joint density of the form
\begin{eqnarray}
   { {h}(x,y)=\left\{\begin{array}{cc}
                    \alpha^{\prime}\beta \exp[{-(\alpha+\beta-\alpha^{\prime})y-\alpha^{\prime}x}], &
                    x\ge
                    y> 0\vspace{0.1cm}\\\alpha\beta^{\prime}\exp[{-(\alpha+\beta-\beta^{\prime})x-\beta^{\prime}y}], & y> x> 0,
                  \end{array}\right.}
\end{eqnarray}
where $\alpha,\alpha^{\prime},\beta,\beta^{\prime}>0.$ If
$\alpha+\beta> \alpha^{\prime}\vee\beta^{\prime}$, its survival
function is equal to
\begin{eqnarray*}
   {\overline{H}(x,y)=\left\{\begin{array}{cc}
                    \frac{\beta}{\alpha+\beta-\alpha^{\prime}}\exp[{-(\alpha+\beta-\alpha^{\prime})y-\alpha^{\prime}x}]+
                    \frac{\alpha-\alpha^{\prime}}{\alpha+\beta-\alpha^{\prime}}\exp[{-(\alpha+\beta)x}], & x\ge
                    y\ge 0\vspace{0.2cm}\\
                    \frac{\alpha}{\alpha+\beta-\beta^{\prime}}
                    \exp[{-(\alpha+\beta-\beta^{\prime})x-\beta^{\prime}y}]+
                    \frac{\beta-\beta^{\prime}}{\alpha+\beta-\beta^{\prime}}\exp[{-(\alpha+\beta)y}], & y\ge x\ge
                    0.
                  \end{array}\right.}
\end{eqnarray*}

  It worths noting that by choosing
  $\alpha=\frac{\lambda_1\lambda}{\lambda_1+\lambda_2},\
  \beta=\frac{\lambda_{2}\lambda}{\lambda_1+\lambda_2},\
  \alpha^{\prime}=\lambda_1+\lambda_{12}$ and
  $\beta^{\prime}=\lambda_2+\lambda_{12},$
Freund's BVE (9) reduces to Block and Basu's BVE (8).
\medskip\\
\noindent{\bf 3. New General Properties of BLM Distributions}

Let $(X,Y)\sim H=BLM(F,G,\theta)\in{\cal BLM}$ with marginals $F$
and $G$ on $(0,\infty)$, parameter $\theta >0$ and survival function
(5). Denote the Laplace-Stieltjes transform of $X$ ($Y,$\,resp.) by
$L_X$ ($L_Y,$\,resp.), and that of $(X,Y)$ by ${\cal L}.$ Then we
have
\smallskip\\
\noindent{\bf Theorem 2.} The {Laplace-Stieltjes transform} of
$(X,Y)\sim H=BLM(F,G,\theta)\in{\cal BLM}$ is
\begin{eqnarray*}{\cal L}(s,t)\equiv
E\left[e^{-sX-tY}\right]=\frac{1}{\theta+s+t}\left[(\theta+s)L_X(s)+(\theta+t)L_Y(t)\right]-\frac{\theta}{\theta+s+t},\
s, t>0.
\end{eqnarray*}

To prove Theorem 2, we need the following lemma due to Lin et
al.\,(2016).
\smallskip\\
 \noindent{\bf
Lemma 1.} Let $(X,Y)\sim H$ defined on ${\mathbb
R}_{+}^{2}=[0,\infty)\times [0,\infty)$. Then the Laplace-Stieltjes
transform of $(X,Y)$\ is equal to\
 $${\cal L}(s,t)=st\int_{0}^{\infty
}\!\!\int_{0}^{\infty }\overline{H} (x,y)e^{-sx-ty}dxdy-1
 +L_X(s)+L_Y(t),~ s,\,t\ge 0.$$
\noindent{\bf Proof of Theorem 2.} We have to calculate the double
integral
\begin{eqnarray*}\int_{0}^{\infty
}\!\!\int_{0}^{\infty }\overline{H}
(x,y)e^{-sx-ty}dxdy=\int\!\!\int_{x\ge y}+\int\!\!\int_{y\ge
x}\equiv A_1+A_2,\end{eqnarray*} where, by changing variables and by integration by parts,
\begin{eqnarray*}
A_1&=&\int_0^{\infty}\!\!e^{-(\theta+t)y}\int_{y}^{\infty}e^{-sx}\overline{F}(x-y)
dxdy
=\int_0^{\infty}\!\!e^{-(\theta+s+t)y}\int_{0}^{\infty}e^{-sz}\overline{F}(z)
dzdy\\
&=&\frac{1}{\theta+s+t}\left[-\frac{1}{s}\int_0^{\infty}\overline{F}(z)de^{-sz}\right]=\frac{1}{\theta+s+t}\left[\frac{1}{s}(1-L_X(s))\right],
\end{eqnarray*}
and similarly,
\begin{eqnarray*}
A_2=\frac{1}{\theta+s+t}\left[\frac{1}{t}(1-L_Y(t))\right].
\end{eqnarray*}
Lemma 1 together with the above $A_1$ and $A_2$ completes the
proof.

 Denote the moment generating function (mgf) of $X$ ($Y,$\,resp.)  by
$M_X$ ($M_Y,$\,resp.), and that of $(X, Y)$ by ${\cal M}.$ Then we
have the following general result.
\smallskip\\
 \noindent{\bf Theorem 3.} Let $(X,Y)\sim
H=BLM(F,G,\theta)\in{\cal BLM}$ and let $r,s$ be real numbers such
that $s+t<\theta.$ Then the {mgf} of $(X,Y)$ is
\begin{eqnarray*}{\cal
M}(s,t)\equiv E\left[e^{sX+tY}\right]
=\frac{1}{\theta-s-t}\left[(\theta-s)M_X(s)+(\theta-t)M_Y(t)\right]-\frac{\theta}{\theta-s-t},
\end{eqnarray*}
provided the expectations (mgfs) exist.

To prove Theorem 3, we need instead the following lemma  due to Lin
et al.\,(2014).
\smallskip\\
\noindent{\bf Lemma 2.} Let $(X,Y)\sim H$ defined on ${\mathbb
R}_{+}^{2}.$
 Let $\alpha $ and $\beta $ be two
increasing and left-continuous functions on ${\mathbb R}_{+}$. Then
the expectation of the product $\alpha (X)\beta (Y)$ is equal to
\begin{eqnarray*}
E[\alpha (X)\beta (Y)]=\int_{0}^{\infty }\!\!\int_{0}^{\infty
}\overline{H} (x,y)d\alpha (x)d\beta (y) -\alpha (0)\beta (0)
+\alpha (0)E[\beta (Y)]+\beta (0)E[\alpha (X)],
\end{eqnarray*}
provided the expectations exist.
\smallskip\\
\noindent{\bf Proof of Theorem 3.} Case (i): $s,t\ge 0.$ Let
$\alpha(x)=e^{sx}$ and $\beta(y)=e^{t y}$ in Lemma 2, then
$${\cal M}(s,t)=st\int_{0}^{\infty }\!\!\int_{0}^{\infty
}\overline{H} (x,y)e^{sx+ty}dxdy-1+M_X(s)+M_Y(t).$$ We have to
calculate  the double integral
\begin{eqnarray*}\int_{0}^{\infty
}\!\!\int_{0}^{\infty }\overline{H}
(x,y)e^{sx+ty}dxdy=\int\!\!\int_{x\ge y}+\int\!\!\int_{y\ge
x}\equiv B_1+B_2,\end{eqnarray*} where, as before,
\begin{eqnarray*}
B_1=\frac{1}{\theta-s-t}\left[\frac{1}{s}(-1+M_X(s))\right]
~~\hbox{and}~~
B_2=\frac{1}{\theta-s-t}\left[\frac{1}{t}(-1+M_Y(t))\right].
\end{eqnarray*}
Lemma 2 together with the above $B_1$ and $B_2$ completes the proof of Case (i).\\
Case (ii): $s\ge 0,\ t< 0.$ To apply Lemma 2, set $\alpha(x)=e^{sx}$
and $\beta(y)=1-e^{ty}.$ Then both $\alpha$ and $\beta$ are
increasing functions on ${\mathbb R}_+$ and $E[\alpha (X)\beta
(Y)]=M_X(s)-{\cal M}(s,t).$ Therefore,
\begin{eqnarray*}{\cal M}(s,t)=M_X(s)-E[\alpha (X)\beta (Y)]
= M_X(s)-1+M_Y(t)+st\int_{0}^{\infty }\!\!\int_{0}^{\infty
}\overline{H} (x,y)e^{sx+ty}dxdy.
\end{eqnarray*}
As before, we carry out the above double integral and complete the
proof of Case (ii).\\ Case (iii): $s< 0,\ t\ge 0.$ Set
$\alpha(x)=1-e^{sx}$ and $\beta(y)=e^{ty}$ in
Lemma 2. The remaining proof is similar to that of Case (ii) and is omitted.\\
Case (iv): $s, t<0.$ This case was treated in Theorem 1. The proof
is completed.

Next, we consider the product moments of BLM distributions.
\smallskip\\
\noindent{\bf Theorem 4.} For positive integers $i$ and $j$, the
{product moment} $E[X^iY^j]$ of $(X,Y)\sim H=BLM(F,G,\theta)\in{\cal
BLM}$ is of the form
\begin{eqnarray*}E[X^iY^j]=i\,j\,\sum_{k=0}^{i-1}\frac{1}{i-k}{i-1\choose
k}\frac{\Gamma(j+k)}{\theta^{j+k}}E[X^{i-k}]
+\,i\,j\,\sum_{k=0}^{j-1}\frac{1}{j-k}{j-1\choose
k}\frac{\Gamma(i+k)}{\theta^{i+k}}E[Y^{j-k}],
\end{eqnarray*}
provided the expectations exist.

The first product moment has a neat representation in terms of
marginal means and the parameter $\theta$, from which we can
calculate Pearson's correlation of BLM distributions.
\smallskip\\
\noindent{\bf Corollary 1.} $E[XY]=\frac{1}{\theta}(E[X]+E[Y])$
provided the expectations exist.

To prove Theorem 4 above, we will apply the following lemma due to
Lin et al.\,(2014).
\smallskip\\
\noindent {\bf Lemma 3.} Let $(X,Y)\sim H$ defined on ${\mathbb
R}_+^2$, and let the expectations $E[X^rY^s]$, $E[X^r]$ and $
E[Y^s]$ be finite for some positive real numbers $r$ and $s$. Then
the product moment
\[
E[X^rY^s]=rs\int_0^{\infty}\!\!\int_0^{\infty}\overline{H}
(x,y)x^{r-1}y^{s-1}dxdy.
\]

\noindent{\bf Proof of Theorem 4.} We have to calculate the double
integral
\begin{eqnarray*}\int_{0}^{\infty
}\!\!\int_{0}^{\infty }\overline{H}
(x,y)x^{i-1}y^{j-1}dxdy=\int\!\!\int_{x\ge y}+\int\!\!\int_{y\ge
x}\equiv C_1+C_2,\end{eqnarray*} where, by changing variables and by integration by parts,
\begin{eqnarray*}
C_1&=&\int_0^{\infty}\!\!e^{-\theta
y}y^{j-1}\int_{y}^{\infty}x^{i-1}\overline{F}(x-y) dxdy
=\int_0^{\infty}\!\!e^{-\theta
y}y^{j-1}\int_{0}^{\infty}(y+z)^{i-1}\overline{F}(z)
dzdy\\
&=&\sum_{k=0}^{i-1}{i-1\choose k}
\int_0^{\infty}\!\!y^{j-1+k}e^{-\theta
y}\int_{0}^{\infty}z^{i-1-k}\overline{F}(z)
dzdy\\
&=&\sum_{k=0}^{i-1}{i-1\choose
k}\frac{\Gamma(j+k)}{\theta^{j+k}}\left[\frac{1}{i-k}\int_0^{\infty}\overline{F}(z)dz^{i-k}\right]=\sum_{k=0}^{i-1}\frac{1}{i-k}{i-1\choose
k}\frac{\Gamma(j+k)}{\theta^{j+k}}E[X^{i-k}],
\end{eqnarray*}
and similarly, $$C_2=\sum_{k=0}^{j-1}\frac{1}{j-k}{j-1\choose
k}\frac{\Gamma(i+k)}{\theta^{i+k}}E[Y^{j-k}].$$ Finally, Lemma 3
together with the above $C_1$ and $C_2$ completes the proof.

For moment generating functions of some {\it specific} BLM
distributions, see Chapter 47 of Kotz et al.\,(2000), while for
product moments of such distributions, see Nadarajah (2006).

For the next and later results, we need some notations in
reliability theory. For random variables $X\sim F$ and $Y\sim G,$ we
say that $X$ is smaller than $Y$ in the usual stochastic order
(denoted by $X\le_{st} Y$) if $\overline{F}(x)\le \overline{G}(x)$
for all $x$, that  $X$ is smaller than $Y$ in the hazard rate order
(denoted by $X\le_{hr} Y$) if $\overline{G}(x)/\overline{F}(x)$ is
increasing in $x$, and that $X$ is smaller than $Y$ in the reversed
hazard rate order (denoted by $X\le_{rh} Y$) if ${G}(x)/{F}(x)$ is
increasing in $x$. Suppose $F$ and $G$ have densities $f$ and $g$,
respectively. Then we say that $X$ is smaller than $Y$ in the
likelihood ratio order (denoted by $X\le_{\ell r} Y$) if $g(x)/f(x)$
is increasing in $x$. For more definitions of the related stochastic
orders, see, e.g., M\"uller and Stoyan (2002), Shaked and
Shanthikumar (2007), Lai and Xie (2006) as well as Kayid et al.
(2016). The latter studied stochastic comparisons of the age
replacement models.

On the other hand, for a distribution $F$ itself we define the
notions of increasing failure rate (IFR), decreasing failure rate
(DFR), increasing failure rate in average (IFRA), and decreasing
failure rate in average (DFRA) as follows. We say that\\ (a) $F$ is
IFR (DFR,\,resp.) if $-\log\overline{F}(x)$ is convex
(concave,\,resp.) in $x\ge 0,$ and\\  (b)  $F$ is IFRA
(DFRA,\,resp.) if $-(1/x)\log \overline{F}(x)$
 is increasing (decreasing,\,resp.) in $x> 0$,
or, equivalently, $\overline{F}^{\alpha}(x)\le
(\ge,\,\hbox{resp.})\,\, \overline{F}(\alpha x)$ for all
$\alpha\in(0,1)$ and $x\ge 0.$ (See Barlow and Proschan 1981,
Chapters 3 and 4.)

The
 bivariate IFRA and bivariate DFRA distributions $H$ can be defined
similarly:\\ $H$ is bivariate IFRA (DFRA,\,resp.) if
$\overline{H}^{\alpha}(x,y)\le (\ge,\,\hbox{resp.})\,\,
\overline{H}(\alpha x, \alpha y)$ for all $\alpha\in(0,1)$ and $x,\
y\ge 0$ (see Block and Savits 1976, 1980). It worths mentioning that
there are some other definitions of bivariate IFRA distributions
that all extend the univariate case (see, e.g., Esary and Marshall
1979 or Shaked and Shanthikumar 1988).

Using reliability language, we have the following useful results.
Especially, Theorem 5(iii) means that in the ${\cal BLM}$ family,
positive bivariate aging plays in favor of positive univariate aging
in the sense of IFRA, and vice versa. This is  in general not true
even under the condition of positive dependence for lifetimes; see
Bassan and Spizzichino (2005, Remark 6.8), which analyzed the
relations among univariate and bivariate agings and dependence.
\smallskip\\
 \noindent{\bf Theorem 5.} Let $(X, Y)\sim
H=BLM(F,G, \theta)\in{\cal BLM}$ and  ${Z\sim Exp(\theta)}$. Then\\
(i) $Z\le_{\ell r} X$ and $Z\le_{\ell r} Y$;\\ (ii) $Z\le_{st} X$
and $Z\le_{st} Y$;
{$Z\le_{hr} X$} and $Z\le_{hr} Y$; {$Z\le_{rh} X$} and $Z\le_{rh} Y$;\\
(iii) $(X,Y)$ has a bivariate  IFRA distribution iff both
marginals $F$ and $G$ are
IFRA;\\
(iv)  $(X,Y)$ has a bivariate  DFRA distribution iff both marginals
$F$ and $G$ are DFRA.\\ {\bf Proof.} Part (i) follows immediately
from Theorem 1(iii) (see Ghurye and Marshall 1984), while part (ii)
follows from the fact that the likelihood ratio order is stronger
than the {usual stochastic order, hazard rate order and reversed
hazard rate order} (M\"uller and Stoyan 2002,\,pp.\,12--13). Part
(iii) holds true by verifying that $\overline{H}(\alpha x,\alpha
y)\ge\overline{H}^{\alpha}(x,y)\ \forall \alpha\in(0,1),\,x,y$ $\ge
0$,\ if, and only if,\ (a) $\overline{F}(\alpha x)\ge
\overline{F}^{\alpha}(x)\ \forall \alpha\in(0,1),\,x\ge 0,$ and (b)
$\overline{G}(\alpha y)\ge \overline{G}^{\alpha}(y)\ \forall
\alpha\in(0,1),\,y\ge 0.$ The proof of part (iv) is similar.

Applying the above stochastic inequalities, we can simplify the
proof of some previous known results. For example, we have
\smallskip\\
\noindent{\bf Corollary 2.}  Let $(X, Y)\sim H=BLM(F,G, \theta)\in{\cal BLM}$. Then the following statements are true.\\
(i) Both hazard rates of marginals $F, G$  are bounded by $\theta$ and hence $\theta\ge f(0)\vee g(0)$.\\
 (ii) Both
the functions $F(-\frac{1}{\theta}\log (1-t))$ and
$G(-\frac{1}{\theta}\log (1-t))$ are convex in $t\in [0,1)$, and
hence
 $f^{\prime}(x)+\theta f(x)\ge 0,\ g^{\prime}(x)+\theta g(x)\ge 0,\ x\ge 0,$\ \ if $f$ and $g$ are differentiable.\\
 (iii) Let
$S_F, S_G$ be the supports of marginals $F, G$ with densities $f,
g$, respectively. Then $S_F=[a_F,\infty), S_G=[a_G,\infty)$ for some
nonnegative constants $a_F, a_G$ with $a_Fa_G=0,$ and $f, g$ are
positive on $(a_F,\infty), (a_G,\infty)$, respectively.\\
(iv) If $H$ is not absolutely continuous, then $f(0)>0$, $g(0)>0,$ and hence $a_F=a_G=0$.\\
{\bf Proof.} Part (i) follows from the facts {$Z\le_{hr} X$} and
$Z\le_{hr} Y$, where ${Z\sim Exp(\theta)}$, while part (ii) is due
to the probability-probability plot characterization for $Z\le_{\ell
r} X$ and $Z\le_{\ell r} Y$ (see Theorem 1.4.3 of M\"uller and
Stoyan 2002).  Part (iii) follows from the facts $Z\le_{\ell r} X$,
$Z\le_{\ell r} Y$ and Theorem 1(iv), because the latter implies that
at least one of the left extremities $a_F$ and $a_G$ of marginal
distributions should be zero. Finally, to prove part (iv), we note
that $\Pr(X-Y>0)=1-f(0)/\theta$ and $\Pr(Y-X>0)=1-g(0)/\theta$ by
Theorem 1(ii) (see Ghurye and Marshall 1984, p.\,789). So if $H$ is
not absolutely continuous, $\Pr(X=Y)>0,$ and hence
$f(0)=\theta\Pr(X\le Y)>0$ and $g(0)=\theta\Pr(Y\le X)>0$. The proof
is complete.
\medskip\\
\noindent{\bf 4. Dependence Structures of BLM Distributions}

 Recall that a bivariate distribution $H$ with marginals $F$ and $G$ is
positively
 quadrant dependent (PQD) if
 $$H(x,y)\ge F(x)G(y)\ \ \forall\, x,y\ge 0,\ \ \hbox{or,
 equivalently,}\ \ \overline{H}(x,y)\ge \overline{F}(x)\overline{G}(y)\ \ \forall\, x,y\ge
 0,$$
 which implies that  $H$ has a nonnegative covariance by Hoeffding
 representation for covariance (see, e.g., Lin et al.\,2014, p.\,2).
 A stronger (positive dependence) property than the PQD
is the total positivity defined below. For a nonnegative function
$K$ on the rectangle $(a,b)\times (c,d)$ (or on the product of two
subsets of ${\mathbb R}$), we say that $K(x,y)$ is {totally positive
of order $r$ (TP$_r,\, r\ge 2$)} in $x$ and $y$ if for each fixed
$s\in\{2,3,\ldots,r\}$ and for all $a<x_1<x_2<\cdots<x_s<b\
\mbox{and} \ c<y_1<y_2<\cdots<y_s<d$, the determinant of the
$s\times s$ matrix $(K(x_i,y_j))$ is nonnegative. The function $K$
is said to be {TP$_{\infty}$} if it is TP$_r$ for any order $r\ge 2$
(Karlin 1968).

The total positivity plays an important role on various concepts of
bivariate dependence (see, e.g., Shaked 1977 and Lee 1985).
Moreover, applying total positivity of the bivariate distribution or
its survival function, we can derive some useful probability
inequalities, among many applications to applied fields including
statistics, reliability and economics (see, e.g.,  Gross and
Richards 1998, 2004, and Karlin and Proschan 1960). Especially, the
latter studied the totally positive kernels that arise from
convolutions of P\'olya type distributions.

 We now characterize the
{TP$_{2}$} property of the survival functions of BLM
distributions.
\smallskip\\
{\bf Theorem 6.} Let $(X,Y)\sim H=BLM(F,G,\theta)\in{\cal BLM}.$
Then the survival function $\overline{H}$ is {TP$_{2}$}\ iff
 the marginal distributions $F$ and $G$ are
IFR and together satisfy $\overline{F}(x)\overline{G}(x)\le
\exp(-\theta
x),\, x\ge 0.$\\
{\bf Proof.} We define the cross-product ratio of $\overline{H}$:
$$r\equiv r(x_1,x_2;y_1,y_2)=\frac{\overline{H}(x_1,y_1)\overline{H}(x_2,y_2)}
{\overline{H}(x_1,y_2)\overline{H}(x_2,y_1)},\ \ 0<x_1<x_2,\
0<y_1<y_2.$$ Then, by definition, $\overline{H}$ is {TP$_{2}$} iff
$r(x_1,x_2;y_1,y_2)\ge 1$ for all $0<x_1<x_2,\ 0<y_1<y_2$.\\
(Necessity) Suppose that $\overline{H}$ is TP$_{2}$. Then for all
$0<x_1=y_1<x_2=y_2$, we have
\[r=
r(x_1,x_2;x_1,x_2)=
\frac{\overline{H}(x_1,x_1)\overline{H}(x_2,x_2)}
{\overline{H}(x_1,x_2)\overline{H}(x_2,x_1)}=\frac{\exp(-\theta
(x_2-x_1)) }{\overline{F}(x_2-x_1)\overline{G}(x_2-x_1)}\ge 1.\]

This implies that $\overline{F}(x)\overline{G}(x)\le \exp(-\theta
x),\, x\ge 0.$ Next, we prove that the marginal distribution $G$ is
IFR. Note that the following statements are equivalent:\\
(i)  $g(y)/\overline{G}(y)$ is increasing in $y\ge 0$,\\
(ii) $\frac{\overline{G}(y+t)}{\overline{G}(t)}$ is decreasing in
$t\in(0,\infty)$ for each $y\ge 0$ (Barlow and
Proschan 1981, p.\,54),\\
(iii) $\frac{\overline{G}(t)}{\overline{G}(y+t)}$
is increasing in $t\in(0,\infty)$ for each $y\ge 0$,\\
(iv) $\frac{\overline{G}(y-x_2)}{\overline{G}(y-x_1)}$ is increasing
in
$y> x_2$ for any fixed $0<x_1<x_2$,\\
(v) the ratio $r^*_{\overline{G}}\equiv
\frac{\overline{G}(y_1-x_1)\overline{G}(y_2-x_2)}
{\overline{G}(y_2-x_1)\overline{G}(y_1-x_2)}\ge 1\ \ \hbox{for all}\
\ 0<x_1<x_2< y_1<y_2.$\\ The latter is true because in this case
$r^*_{\overline{G}}=r(x_1,x_2;y_1,y_2)\ge 1$ by (5) and the
assumption. Similarly, we can prove that $F$ is IFR because the
ratio
$$r^*_{\overline{F}}\equiv\frac{\overline{F}(x_1-y_1)\overline{F}(x_2-y_2)}
{\overline{F}(x_2-y_1)\overline{F}(x_1-y_2)}\ge 1\ \ \hbox{for all}\
\ 0<y_1<y_2< x_1<x_2.$$

(Sufficiency) Suppose that the marginal distributions $F$ and $G$
are IFR and together satisfy $\overline{F}(x)\overline{G}(x)\le
\exp(-\theta x),\, x\ge 0.$ Then we want to prove that
$\overline{H}$ is TP$_{2}$, that is, for all $0<x_1<x_2,\ 0<y_1<
y_2$, the cross-product ratio $r=r(x_1,x_2;y_1,y_2)\ge 1.$ Without
loss of generality, we consider only three possible cases below,
\[(a)\ 0<x_1\le x_2\le y_1\le y_2,\ \ (b)\ 0<x_1\le y_1\le x_2\le
y_2,\ \ (c)\ 0<x_1\le y_1\le y_2\le x_2,\] because the remaining
cases can be proved by exchanging the roles of $F$ and $G$.

For case (a), we have $r\ge 1$ by the equivalence relations shown in
the necessity part and by the continuity of $H$ when $x_2=y_1$. For
case (b),  the cross-product ratio
\[r=\frac{\exp(-\theta x_2)\overline{G}(y_1-x_1)\overline{G}(y_2-x_2)}
{\exp(-\theta y_1)\overline{G}(y_2-x_1)\overline{F}(x_2-y_1)}\ge
\frac{\overline{G}(y_1-x_1)\overline{G}(y_2-x_2)\overline{G}(x_2-y_1)}
{\overline{G}(y_2-x_1)},\] because
$\overline{F}(x_2-y_1)\overline{G}(x_2-y_1)\le \exp(-\theta
(x_2-y_1))$ by the assumption. Recall that any IFR distribution is
new better than used (Barlow and Proschan 1981, p.\,159).
Therefore, $\overline{G}(x+y)\le \overline{G}(x)\overline{G}(y)$
for all $x,y\ge 0,$ and hence the last $r\ge 1.$ Similarly, for
case (c),
\[r=\frac{\exp(-\theta y_2)\overline{G}(y_1-x_1)\overline{F}(x_2-y_2)}
{\exp(-\theta y_1)\overline{G}(y_2-x_1)\overline{F}(x_2-y_1)}\ge
\frac{\overline{G}(y_1-x_1)\overline{G}(y_2-y_1)}{\overline{G}(y_2-x_1)}
\times\frac{\overline{F}(y_2-y_1)\overline{F}(x_2-y_2)}{\overline{F}(x_2-y_1)}\ge
1,
\]
by the assumptions. This completes the proof.

Recall also that
 for any bivariate distribution $H$ with
marginals $F$ and $G$, there exist a {copula $C$} (a bivariate
distribution with uniform marginals on $[0,1]$) and a {survival
copula $\hat{C}$} such that $H(x,y)=C(F(x), G(y))$ and
$\overline{H}(x,y)=\hat{C}(\overline{F}(x),\overline{G}(y))$ for all
$x,y\in {\mathbb R}\equiv (-\infty,\infty).$ Namely, $C$ links $H$
and $(F,G)$, while $\hat{C}$ links $\overline{H}$ and
$(\overline{F},\overline{G})$.
\smallskip\\
\noindent{\bf Corollary 3.}  Let $(X,Y)\sim
H=BLM(F,G,\theta)\in{\cal BLM}.$ Then
 the survival copula $\hat{C}$ of $H$ is {TP$_{2}$}\ iff
 the marginal distributions $F$ and $G$ are
IFR and together satisfy $\overline{F}(x)\overline{G}(x)\le
\exp(-\theta
x),\, x\ge 0.$\\
{\bf Proof.} Since the marginal $F$ is absolutely continuous on the
support $[a_F,\infty)$ with positive density $f$ on $(a_F,\infty)$
(see Corollary 2(iii) above), $F$ is strictly increasing and
continuous on $(a_F,\infty)$. Similarly,  the marginal $G$  is
strictly increasing and continuous on $(a_G,\infty)$. By Theorem 6,
it suffices  to prove that $\overline{H}$ is {TP$_{2}$} on
$(a_F,\infty)\times (a_G,\infty)$ iff its survival copula $\hat{C}$
is  {TP$_{2}$} on $(0,1)^2$. Recall the facts (i)
$\overline{H}(x,y)=\hat{C}(\overline{F}(x),\overline{G}(y))$,
$(x,y)\in (a_F,\infty)\times (a_G,\infty)$, (ii)
$\hat{C}(u,v)=\overline{H}(\overline{F}^{\,-1}(u),
\overline{G}^{\,-1}(v))$, $u,v\in(0,1)$, where $\overline{F}^{\,-1},
\overline{G}^{\,-1}$ are inverse functions of $\overline{F},
\overline{G}$, respectively, and (iii) all the functions
$\overline{F}, \overline{G}, \overline{F}^{\,-1}$ and
$\overline{G}^{\,-1}$ are decreasing. The required result then
follows immediately (see, e.g., Lemma 5(ii) below).

The counterpart of TP$_{2}$ property is the reverse regular of order
two (RR$_2$). For a nonnegative function $K$ on $(a,b)\times(c,d)$,
we say that $K$  is RR$_2$ if the determinant of the $2\times 2$
matrix $(K(x_i,y_j))$ is non-positive for all $a<x_1<x_2<b$ and
$c<y_1<y_2<d$ (see, e.g., Esna--Ashari and Asadi 2016 for examples
of  RR$_2$ joint  densities and survival functions). Mimicking the
proof of Theorem 6, we conclude that for $H=BLM(F,G,\theta)\in{\cal
BLM}$, the survival function $\overline{H}$ is RR$_2$\ iff  the
survival copula $\hat{C}$ of $H$ is RR$_2$\ iff
 the marginal distributions $F$ and $G$ are
DFR and satisfy $\overline{F}(x)\overline{G}(x)\ge \exp(-\theta
x),\, x\ge 0.$ To construct such a BLM distribution with RR$_2$
survival function, we first consider the Pareto Type II distribution
(or Lomax distribution) $F$  with density function
$f(x)=(\alpha/\beta) (1+x/\beta)^{-(\alpha+1)},\ x\ge 0,$ and
survival function $\overline{F}(x)=(1+x/\beta)^{-\alpha},\ x\ge 0,$
where $\alpha,\ \beta>0.$ Then choose the parameters: $\alpha\ge 1,\
\beta>0$ and $\theta=(\alpha+1)/\beta.$ It can be checked that the
$\overline{H}$ defined in (5) with $G=F$ is a {\it bona fide}
survival function, and is RR$_2$ if $\alpha=1.$

It is seen that all the conditions in Theorem 6 are satisfied by the
Marshall--Olkin BVE. Therefore, the survival function and survival
copula of the Marshall--Olkin BVE are both TP$_{2}$, regardless of
parameters; a more general result will be given in Theorem 8 below.
We next characterize, by a different approach, the TP$_{2}$ property
of some joint densities of absolutely continuous BLM distributions.
\smallskip\\
\noindent {\bf Theorem 7.} Let $H=BLM(F,G,\theta)\in{\cal BLM}$ be
absolutely continuous and have joint density function $h$. Suppose
that the marginal density functions $f$ and $g$ are three times
differentiable on $(0,\infty)$ and  that $\theta
f(0^+)+f^{\prime}(0^+)=\theta g(0^+)+g^{\prime}(0^+)$ is finite.
Assume further the functions
$$h_1(x|\theta)\equiv\theta f(x)+f^{\prime}(x)>0,\ \ x>0,\ \ \hbox{and}\ \
h_2(y|\theta)\equiv\theta g(y)+g^{\prime}(y)>0,\ \ y>0.$$ Then the
joint density function $h$ is {TP$_{2}$}\ iff the marginal densities
satisfy (i) $\left(h_i^{\prime}(x|\theta)\right)^2\ge
h_i^{\prime\prime}(x|\theta)h_i(x|\theta),\, x>0,\,i=1,2,$ and (ii)
$h_1(x|\theta)h_2(x|\theta)\le h_1^2(0^+|\theta)\exp(-\theta x),\,
x>0$.

 To prove this theorem, we need the concept of local dependence
function and the following lemma, in which part (ii) is essentially
due to Holland and Wang (1987, p.\,872). An alternative (complete)
proof of part (ii) is provided below. In their proof, Holland and
Wang (1987) assumed implicitly the {\it integrability} of the local
dependence function, while Kemperman (1977, p.\,329) gave without
proof the same result under {\it continuity} (smoothness) condition
(see also Newman 1984). Wang (1993) proved that a positive
continuous bivariate density on a Cartesian product
$(a,b)\times(c,d)$ is uniquely determined by its marginal densities
and local dependence function when the latter exists and is
integrable. On the other hand, Jones (1996, 1998) investigated the
bivariate distributions with constant local dependence.
\smallskip\\
\noindent{\bf Lemma 4.} Let $K$ be a positive function on
$D=(a,b)\times(c,d)$. Then we have\\
(i) $K$ is {TP$_{2}$} on $D$ iff $\log K$ is 2-increasing;\\
(ii) $K$ is {TP$_{2}$} on $D$ iff the local dependence function
$\gamma_K(x,y)\equiv \frac{\partial^2}{\partial x\partial y}\log
K(x,y)\ge 0~~\hbox{on}\ D,$ provided the second-order
partial derivatives exist.\\
{\bf Proof.}  Part (i) is trivial by the definition of 2-increasing
functions (see Nelsen 2006,\,p.\,8), and part (ii) follows from part
(i) and the fact that  under the smoothness assumption, $\log K$ is
2-increasing iff the local dependence function $\gamma_K(x,y)\ge 0.$
To prove part (ii) directly, note that the following statements are
equivalent: (a) $\frac{\partial^2}{\partial x\partial y}\log
K(x,y)\ge 0$ on $D$,\ (b) $\frac{\partial}{\partial y}\log
[K(x_2,y)/K(x_1,y)]\ge 0$ for all $y$ and for all $x_1<x_2$,\ (c)
$\log [K(x_2,y)/K(x_1,y)]$ is increasing in $y$ for all $x_1<x_2$,\
(d) $K(x_2,y)/K(x_1,y)$ is increasing in $y$ for all $x_1<x_2$,\ (e)
$K(x_2,y_2)/K(x_1,y_2)\ge K(x_2,y_1)/K(x_1,y_1)$ for all $y_1<y_2,\
x_1<x_2$,\  (f) the cross-product ratio of $K$ satisfies:
$K(x_1,y_1)K(x_2,y_2)/[K(x_1,y_2)K(x_2,y_1)]\ge 1$ for all
$x_1<x_2,\ y_1<y_2$,\ and (g) the function $K$ is {TP$_{2}$} on $D$.
The proof is complete.
\smallskip\\
{\bf Proof of Theorem 7.}  By the
assumptions, the joint density function of $H$ is of the form
\begin{eqnarray*}
    h(x,y)=\left\{\begin{array}{cc}
                    e^{-\theta y}\,h_1(x-y|\theta), & x\ge y \vspace{0.1cm}\\
                    e^{-\theta x}\,h_2(y-x|\theta),  & x\le y,
                  \end{array}\right.
\end{eqnarray*}
where $h_i(0|\theta)\equiv h_i(0^+|\theta),\,i=1,2.$  For $x\ne y$,
the local dependence function of $h$ is
\begin{eqnarray*}
   \gamma_h(x,y)=\frac{\partial^2}{\partial x\partial
y}\log h(x,y)=\left\{\begin{array}{cc}
                    \frac{[h_1^{\prime}(x-y|\theta)]^2-h_1^{\prime\prime}(x-y|\theta)h_1(x-y|\theta)}
                    {h_1^2(x-y|\theta)}, & x> y \vspace{0.2cm}\\
                    \frac{[h_2^{\prime}(y-x|\theta)]^2-h_2^{\prime\prime}(y-x|\theta)h_2(y-x|\theta)}
                    {h_2^2(y-x|\theta)},  & x< y.
                  \end{array}\right.
\end{eqnarray*}
Therefore, $\gamma_h(x,y)\ge 0$ for all $(x,y)$ with $x\ne y$ iff
the property (i) holds true.

 (Necessity) If $h$ is {TP$_{2}$} on
$(0,\infty)^2$, then it is also {TP$_{2}$} on each rectangle
(rectangular area) in the region ${\cal A}_1=\{(x,y): x>y>0\}$ or in
${\cal A}_2=\{(x,y): y>x>0\}$, and hence the property (i) holds true
by Lemma 4 and the above observation. Next, the property (ii)
follows from the fact that for all $0<x_1=y_1<x_2=y_2$, the
cross-product ratio $r_h$ of $h$ satisfies
\[1\le r_h\equiv r_h(x_1,x_2;y_1,y_2)=\frac{h(x_1,y_1)h(x_2,y_2)}{h(x_1,y_2)h(x_2,y_1)}
=\frac{\exp(-\theta(x_2-x_1))h_1(0|\theta)h_2(0|\theta)}{h_1(x_2-x_1|\theta)h_2(x_2-x_1|\theta)}.\]
This completes the proof of the necessity part.

 (Sufficiency)
Suppose $0<x_1<x_2$ and $0<y_1<y_2$, then we want to prove the
cross-product ratio $r_h\ge 1$ under the assumptions (i) and (ii).
If the rectangle with four vertices $P_i, i=1,2,3,4$, where
$P_1=(x_1,y_1), P_2=(x_2, y_1), P_3=(x_2,y_2), P_4=(x_1,y_2)$, lies
entirely in the region ${\cal A}_1$ or ${\cal A}_2$, then $r_h\ge 1$
by the assumption (i) and Lemma 4. If $0<x_1=y_1<x_2=y_2$, then the
assumption (ii) implies $r_h\ge 1$. For the remaining cases, we
apply the technique of factorization of the cross-product ratio  if
necessary. For example, if $P_*=(x_1,y_*)\in \overline{P_1P_4}$ and
$P^*=(x^*,y_2)\in \overline{P_4P_3}$ denote the intersection of the
diagonal line $x=y$ and boundary of the rectangle, where
$x_1<x^*<x_2$ and $y_1<y_*<y_2$, then we split the original
rectangle
 into four sub-rectangles by adding the
new point $(x^*,y_*)$ and calculate the ratio
\begin{eqnarray*}
r_h(x_1,x_2;y_1,y_2)=r_h(x_1,x^*;y_1,y_*)r_h(x_1,x^*;y_*,y_2)r_h(x^*,x_2;y_1,y_*)r_h(x^*,x_2;y_*,y_2)\ge
1,
\end{eqnarray*}
each factor being greater than or equal to one by the previous
results. The proof is complete.

It is known that the Marshall--Olkin BVE (6) is {PQD}, so are its
copula $C$ and survival copula $\hat{C}$ (see Barlow and Proschan
1981, p.\,129). Moreover, $\overline H$ and $\hat{C}$ are TP$_2$ due
to Theorem 6 and its corollary (see also Nelsen 2006, p.\,163, for a
direct proof) and both
 are  even TP$_{\infty}$ if
$\lambda_1=\lambda_2$ (Lin et al.\,2016). We are now able to extend
these results to the following.
\smallskip\\
\noindent {\bf Theorem 8.} The Marshall--Olkin {survival function}
$\overline{H}$ and {survival copula} $\hat{C}$  are both
TP$_{\infty},$ regardless of parameters.

To prove this theorem, we need two more useful lemmas. Lemma 5 is
well-known (see, e.g., Marshall et al.\,2011, p.\,758), while Lemma
6 is essentially due to Gantmacher and Krein (2002), pp.\,78--79
(see also Karlin 1968, p.\,112, for an alternative version).
\smallskip\\
\noindent{\bf Lemma 5.} Let $r\ge 2$ be an integer.\\
(i) If $k(x,y)$ is TP$_r$ in $x$ and $y$, and if both $u$ and $v$
are nonnegative functions, then the product function
$K(x,y)=u(x)\,v(y)\,k(x,y)$ is TP$_r$ in $x$ and $y$.\\
(ii) If $k(x,y)$ is TP$_r$ in $x$ and $y$, and if $u$ and $v$ are
both increasing, or both decreasing,  then the composition function
$K(x,y)=k(u(x),v(y))$ is TP$_r$ in $x$ and $y$.
\smallskip\\
\noindent{\bf Lemma 6.} Let $\phi$ and $\psi$ be two positive
functions on $(a,b).$ Define the symmetric function
\begin{eqnarray*}
    {K_s}(x,y)=\left\{\begin{array}{cc}
                          \psi(x)\,\phi(y), &\,\, a<y\le
                          x<b\vspace{0.1cm}\\
                          \phi(x)\,\psi(y), &\,\, a<x\le y<b.
                        \end{array}\right.
\end{eqnarray*}
If $\phi(x)/\psi(x)$ is nondecreasing in $x\in (a,b)$, then the
function $K_s(x,y)$ is TP$_{\infty}$ in $x$ and $y$.
\smallskip\\
\noindent{\bf Proof of Theorem 8.}  We prove first that the
Marshall--Olkin {survival function} $\overline{H}$ is TP$_{\infty}.$
Rewrite the survival function (6) as
\begin{eqnarray*}
    \overline{H}(x,y)&=&\left\{\begin{array}{cc}
                    \exp[-(\lambda_1+\lambda_{12})x-\lambda_{2}y], & x\ge y \vspace{0.1cm}\\
                    \exp[-(\lambda_2+\lambda_{12})y-\lambda_{1}x],  & x\le
                    y
                  \end{array}\right.\vspace{0.1cm}\\
                &=&\exp[-\lambda_1x-\lambda_2y]K_s(x,y),
\end{eqnarray*}
where the symmetric function
\begin{eqnarray}K_s(x,y)=\left\{\begin{array}{cc}
                    \exp(-\lambda_{12}x), & x\ge y \vspace{0.1cm}\\
                    \exp(-\lambda_{12}y),  & x\le
                    y.
                  \end{array}\right.
\end{eqnarray}

Let $\phi(x)=1$ and $\psi(y)=\exp(-\lambda_{12}y)$. Then by Lemma 6,
we see that the function $K_s$ in (10) is TP$_{\infty},$ so is
$\overline{H}$ by Lemma 5(i). Next, recall that the Marshall--Olkin
survival copula
$$\hat{C}(u,v)=\overline{H}(\overline{F}^{\,-1}(u),\ \overline{G}^{\,-1}(v)),\ \ \
u,v\in (0,1),$$ where $\overline{F}^{\,-1}$ and $\overline{G}^{\,-1}$ are the inverse (decreasing) functions of
$\overline{F}(x)=\exp[-(\lambda_1+\lambda_{12})x]$ and
$\overline{G}(y)=\exp[-(\lambda_2+\lambda_{12})y]$, respectively.
Therefore, $\hat{C}$ is TP$_{\infty}$ by Lemma 5(ii).

It is well known that if a bivariate distribution $H$ has TP$_2$
density, then both $H$ and its joint survival function
$\overline{H}$ are TP$_2$ (see, e.g., Balakrishnan and Lai 2009,
p.\,116). A more general result is given as follows.
\smallskip\\
\noindent{\bf Theorem 9.} If the bivariate distribution $H$ has
TP$_r$ density with $r\ge 2$, then both $H$ and $\overline{H}$ are
TP$_r$. Consequently, if $H$ has TP$_{\infty}$ density,
then both $H$ and $\overline{H}$ are TP$_{\infty}$.\\
{\bf Proof.} Let us consider first the TP$_{\infty}$ indicator
functions $K_1(x, y)=\hbox{I}_{(-\infty, x]}(y)$ and
$K_2(x,y)=\hbox{I}_{[x,\infty)}(y)$, and then apply Theorem 3.5 of
Gross and Richards (1998) restated below.

For example, to prove the TP$_r$ property of $H$, we have to claim
that {for all $x_1<\cdots<x_r$ and $y_1<\cdots<y_r$, the determinant
of each $s\times s$ sub-matrix $(H(x_i,y_j))$ (with $2\le s\le r$)
is nonnegative. To prove this, let us recall that
$H(x_i,y_j)=E[\hbox{I}_{(-\infty, x_i]}(X)\,\hbox{I}_{(-\infty,
y_j]}(Y)]=E[\phi(i,X)\,\psi(j,Y)]$, where
$\phi(i,x)=\hbox{I}_{(-\infty, x_i]}(x)$ is TP$_r$ in two variables
$i\in \{1,2,\ldots, r\}$ and $x\in {\mathbb R},$ and
$\psi(j,y)=\hbox{I}_{(-\infty, y_j]}(y)$ is TP$_r$ in two variables
$j\in \{1,2,\ldots, r\}$ and $y\in {\mathbb R}$. Then Gross and
Richards' Theorem applies and hence $H$ is TP$_r$. Similarly, the
survival function $\overline H$ is TP$_r$. The proof is complete.
\smallskip\\
\noindent{\bf Gross and Richards' (1998) Theorem.} Let $r\ge 2$ be
an
 integer and let the bivariate $(X,Y)\sim H$ have TP$_r$ density.
Assume further that both the functions $\phi(i,x)$ and $\psi(i,x)$
are TP$_r$ in two variables $i\in \{1,2,\ldots, r\}$ and $x\in
{\mathbb R}$. Then the $r\times r$ matrix
$(E[\phi(i,X)\,\psi(j,Y)])$ is totally positive, that is, all its
minors (of orders $\le r$) are nonnegative real numbers.

 As mentioned in Balakrishnan and Lai (2009, p.\,124),
the Block--Basu BVE (8) is PQD if $\lambda_1=\lambda_2.$ We now
extend this result to the following.
\smallskip\\
\noindent{\bf Theorem 10.} (i) If {$\lambda_1=\lambda_2$} in (8),
then the
Block--Basu BVE has TP$_{\infty}$ density.\\
(ii)  If $\alpha=\beta\le \alpha^{\prime}=\beta^{\prime}$ in (9),
then
the Freund BVE has TP$_{\infty}$ density.\\
{\bf Proof.} Take $\phi(x)=c_1\exp({-\lambda_1x})$ and
$\psi(y)=c_2\exp[{-(\lambda_2+\lambda_{12})y}]$ for some constants
$c_1, c_2>0.$ Then part (i) follows from (8) and Lemma 6. Part (ii)
can be proved similarly.
\smallskip\\
\noindent {\bf Remark 4.} The same approach applies to other
bivariate (non-BLM) distributions like Li and Pellerey's (2011)
generalized Marshall--Olkin bivariate distribution described below.

In Marshall and Olkin's (1967) {shock model:} $(X, Y)= (X_1\wedge
X_3,\, X_2\wedge X_3)$, we assume instead that $X_1, X_2, X_3$ are
independent general positive random variables (not limited to
exponential ones) and that $X_i\sim F_i,\ i=1,2,3.$ Let $R_i=-\log
\overline{F}{\!}_i$ be the hazard function of $X_i.$ Then the
generalized Marshall--Olkin bivariate distribution $H$ has survival
function
\begin{eqnarray}{\overline{H}(x,y)}&=&\Pr(X>x,\,Y>y)
= \Pr(X_1>x,\,X_2>y,\, X_3>\max\{x,y\})\nonumber\\
&=&\exp [-R_1(x)-R_2(y)-R_3(\max\{x,y\})] ,\ \ \ x,y\ge 0,
\end{eqnarray}
which is PQD (Li and Pellerey 2011). (For other related shock
models, see Marshall and Olkin 1967,  Ghurye and Marshall 1984 as
well as Aven  and Jensen 2013, Section 5.3.4.) We now extend this
result and Theorem 8 as follows.
\smallskip\\
\noindent{\bf Theorem 11.} Let $H$ be the generalized
Marshall--Olkin distribution defined in (11). Then\\
 (i)  the survival function $\overline{H}$
 is
 TP$_{\infty};$\\
 (ii)  the
survival copula $\hat{C}$ of $H$ is TP$_{\infty}$, provided the
functions $\overline{F}{\!}_1\overline{F}{\!}_3$ and
$\overline{F}{\!}_2\overline{F}{\!}_3$
are both strictly decreasing.\\
{\bf Proof.} Write the survival function (11) as
\begin{eqnarray*}
    \overline{H}(x,y)&=&\left\{\begin{array}{cc}
                    \exp[-(R_1(x)+R_3(x))-R_2(y)], & x\ge y \vspace{0.1cm}\\
                    \exp[-(R_2(y)+R_3(y))-R_1(x)],  & x\le y
                  \end{array}\right.\\
                  &=&\exp[-R_1(x)-R_2(y)]K_s(x,y),
\end{eqnarray*}
where the  symmetric function
\begin{eqnarray}K_s(x,y)=\left\{\begin{array}{cc}
                    \exp[-R_3(x)], & x\ge y \vspace{0.1cm}\\
                    \exp[-R_3(y)],  & x\le y.
                  \end{array}\right.
\end{eqnarray}

By taking $\phi(x)=1$ and $\psi(y)=\exp[-R_3(y)]$ in Lemma 6, we
know that the function $K_s$ in (12) is TP$_{\infty}$, and hence the
survival function $\overline{H}$ is TP$_{\infty}$ by Lemma 5(i).
This proves part (i).  To prove part (ii), we note that the marginal
survival functions of $H$ are
$\overline{F}(x)=\exp[-\tilde{R}_1(x)],\ x\ge 0,$\ and
$\overline{G}(y)=\exp[-\tilde{R}_2(y)],\ y\ge 0,$ where the two
functions $\tilde{R}_1(x)=R_1(x)+R_3(x),\ x\ge 0,$\ and
$\tilde{R}_2(y)=R_2(y)+R_3(y),\ y\ge 0,$\ are strictly increasing
 by the conditions on $F_i,\ i=1,2,3.$ This in turn
implies that the marginal distribution functions ${F}$ and ${G}$ are
strictly increasing and hence the survival copula
$$\hat{C}(u,v)=\overline{H}(F^{-1}(1-u), G^{-1}(1-v)),\ \ \ u,v\in
(0,1),$$ because $F^{-1}(F(t))=t,\ t\in (0,1),$ where the quantile
function $F^{-1}(t)=\inf\{x: F(x)\ge t\},\, t\in(0,1)$ (see, e.g.,
Shorack and Wellner 1986, p.\,6). Therefore, $\hat{C}$ is
TP$_{\infty}$ by part (i) and Lemma 5(ii). The proof is complete.
\medskip\\
\noindent{\bf 5.  Stochastic Comparisons of BLM Distributions}

To provide more information about BLM distributions, we can study
{stochastic comparisons in the ${\cal BLM}$ family}. As usual,
define the notions of the upper orthant order ($\le_{uo}$), the
concordance order ($\le_{c}$) and the Laplace transform order
($\le_{Lt}$) as follows. Let $(X_i,Y_i)\sim H_i$ with marginals
$(F_i,G_i)$, $i=1,2,$ on ${\mathbb R}_+.$ Then denote (i)
$(X_1,Y_1)\le_{uo} (X_2,Y_2)$ if $\overline{H}{\!}_1(x,y)\le
\overline{H}{\!}_2(x,y)$\ for all $x,y\ge 0$, (ii)
 $(X_1,Y_1)\le_{c}(X_2,Y_2)$
if $(F_1,G_1)=(F_2,G_2)$ and $(X_1,Y_1)\le_{uo} (X_2,Y_2)$, and
(iii)  $(X_1,Y_1)\le_{Lt}(X_2,Y_2)$ if ${\cal L}_1(s,t)\ge {\cal
L}_2(s,t)$ for all $s, t\ge 0$ (M\"uller and Stoyan 2002, Shaked and
Shanthikumar 2007).  We have, for example, the following results
whose proofs are straightforward and are omitted.
\smallskip\\
\noindent{\bf Theorem 12.} Let $(X_i,Y_i)\sim
H_i=BLM(F_i,G_i,\theta_i)\in{\cal BLM},$ $i=1,2.$ Then we have \\
(i) $X_1\le_{st}X_2,$ $Y_1\le_{st}Y_2$ and $\theta_1\ge \theta_2,$\
iff\ $(X_1,Y_1)\le_{uo}(X_2,Y_2)$,  or, equivalently,
$E[K(X_1,Y_1)]\le
E[K(X_2,Y_2)]$ for any bivariate distribution $K$ on ${\mathbb R}_+^2$;\\
(ii) $F_1=F_2, G_1=G_2$ and $\theta_1\ge \theta_2,$\ iff\
$(X_1,Y_1)\le_{c}(X_2,Y_2)$, or, equivalently,
$E[k_1(X_1)k_2(Y_1)]\le E[k_1(X_2)k_2(Y_2)]$ for all increasing
functions
$k_1,k_2$, provided the expectations exist; and\\
(iii) if $X_1\le_{Lt}X_2,$ $Y_1\le_{Lt}Y_2$ and $\theta_1=\theta_2$,
then $(X_1,Y_1)\le_{Lt}(X_2,Y_2)$, or, equivalently,
$E[k_1(X_1)k_2(Y_1)]\ge E[k_1(X_2)k_2(Y_2)]$ for all completely
monotone functions $k_1,k_2$, provided the expectations exist.

When $H_1$ and $H_2$ have the same pair of marginals $(F,G)$,
Theorem 12(i) reduces, by Corollary 1, to the following interesting
result which is related to the famous {Slepian's inequality} for
bivariate normal distributions (see the discussion in Remark 5
below).
\smallskip\\
\noindent{\bf Corollary 4.} Let $(X_i,Y_i)\sim
H_i=BLM(F,G,\theta_i)\in{\cal BLM},$  with correlation $\rho_i$,
$i=1,2.$  Then {$\rho_1\le\rho_2$} { iff }
{$\overline{H}{\!}_1(x,y)\le \overline{H}{\!}_2(x,y)$}  for all
$x,y\ge 0,$ or, equivalently, {${H}_1(x,y)\le {H}_2(x,y)$}
 for all $x,y\ge 0.$
\smallskip\\
 \noindent{\bf Remark 5.} In Corollary 4 above, if we consider standard bivariate
normal distributions instead of BLM ones, then the conclusion also
holds true and the necessary part is the so-called Slepian's
lemma/inequality; see Slepian (1962), M\"uller and Stoyan (2002),
p.\,97, and Hoffmann-J$\o$rgensen (2013) for more general results.
In Wikipedia, it was said that while this intuitive-seeming result
is true for Gaussian processes,  {it is not in general true for
other random variables}. However, as we can see in Corollary 4,
there are
 infinitely many BLM distributions sharing the same
Slepian's inequality with bivariate normal ones.
\smallskip\\
\noindent{\bf Remark 6.}  We finally compare the effects of the
dependence structure of BLM distributions in different coherent
systems. Consider a two-component system and let the two components
have lifetimes $(X,Y)\sim H=BLM(F,G,\theta)$. Then the lifetime of a
series system composed of these two components is $X\wedge Y\sim
Exp(\theta)$, while the lifetime of a parallel system composed of
the same components is $X\vee Y$ obeying the distribution
$H_p(z)=e^{-\theta z}-1+F(z)+G(z),\ z\ge 0.$ Therefore the mean
times to failure of series and parallel systems are, respectively,
$E[X\wedge Y]=\int_0^{\infty}\overline{H}(x,x)dx=1/\theta$
(decreasing in $\theta$) and $E[X\vee
Y]=E[X]+E[Y]-\int_0^{\infty}\overline{H}(x,x)dx=E[X]+E[Y]-1/\theta$
(increasing in $\theta$).
 The latter further implies that $\theta\ge(E[X]+E[Y])^{-1}$ (compare with Theorem
1(vi)) and that $E[XY]\in[1/\theta^2, (E[X]+E[Y])^2]$ by Corollary
1, provided the expectations exist. See also Aven and Jensen (2013,
Section 2.3) for special cases with exponential marginals as well as
Lai and Lin (2014) for more general results.
\smallskip\\
\noindent{\bf Acknowledgments.} The authors would like to thank the
Editor-in-Chief and two Referees for helpful comments and
constructive suggestions which improve the presentation of the
paper. The paper was presented at the Ibusuki International Seminar,
Ibusuki Phoenix Hotel, held from 6th to 8th March 2016 by Waseda
University, Japan. The authors thank the organizer Professor
Masanobu Taniguchi for his kind invitation and the audiences for
their comments and suggestions.
\bigskip\\
\noindent{\bf References}
\begin{description}

\item Aven, T. and Jensen, U. (2013). {\it Stochastic Models in Reliability}, 2nd ed. Springer, New York.

\item {Balakrishnan, N. and Lai, C.-D. (2009)}. {\it Continuous
Bivariate Distributions}, 2nd ed. Springer, New York.

\item {Barlow, R.\,E. and Proschan, F. (1981)}. {\it
Statistical Theory of Reliability and Life Testing: Probability
Models}, To Begin With. Silver Spring, MD.

\item Bassan, B. and Spizzichino, F. (2005). Relations among
univariate aging, bivariate aging and dependence for
exchangeable lifetimes. {\it J. Multivariate Anal.}, {\bf 93},
313--330.

\item {Block, H.\,W. (1977)}. A characterization of a bivariate
exponential distribution. {\it Ann. Statist.}, {\bf 5}, 808--812.

\item {Block, H.\,W. and Basu, A.\,P. (1974)}. A continuous
bivariate exponential extension. {\em J. Amer. Statist. Assoc.},
{\bf 69}, 1031--1037.

\item Block, H.\,W.  and Savits, T.\,H. (1976). The IFRA closure
problem. {\it Ann. Probab.}, {\bf 4}, 1030--1032.

\item Block, H.\,W.  and Savits, T.\,H. (1980). Multivariate
increasing failure rate average distributions. {\it Ann. Probab.},
{\bf 8}, 793--801.

\item Crawford, G. B. (1966). Characterization of geometric and
exponential distributions. {\it Ann. Math. Statist.}, {\bf 37},
1790--1795.

\item Esary, J.\,D. and Marshall, A.\,W. (1979). Multivariate distributions
with increasing hazard rate average. {\it Ann. Probab.}, {\bf
7}, 359--370.

\item Esna--Ashari, M. and Asadi, M. (2016). On additive--multiplicative
hazards model. {\it Statistics}, {\bf 50}, 1421--1433.

\item {Feller, W. (1965)}. {\it An Introduction to Probability
Theory}, Vol. I. Wiley, New York.

\item {Ferguson, T.\,S. (1964)}. A characterization of the
exponential distribution. {\it Ann. Math. Statist.}, {\bf 35},
1199--1207.

\item {Ferguson, T.\,S. (1965)}. A characterization of the
geometric distribution. {\it Amer. Math. Monthly}, {\bf 72},
256--260.

\item {Fortet, R. (1977)}. {\it Elements of Probability
Theory}.  Gordon and Breach, New York.

\item {Freund, J.\,E. (1961)}. A bivariate extension of the
exponential distribution. {\it J. Amer. Statist. Assoc.}, {\bf
56}, 971--977.

\item Galambos, J. and Kotz, S. (1978). {\it Characterizations of
Probability Distributions.} Springer, New York.

\item {Gantmacher, F.\,R. and Krein, M.\,G. (2002)}. {\it
Oscillation Matrices and Kernels and Small Vibrations of
Mechanical Systems}, Revised edn. Translation based on the 1941
Russian original, Providence, RI.

\item {Ghurye, S.\,G.  (1987)}. Some multivariate lifetime distributions.
{\em Adv. Appl. Probab.}, {\bf 19}, 138--155.

\item {Ghurye, S.\,G. and Marshall, A.\,W. (1984)}. Shock
processes with aftereffects and multivariate lack of memory. {\em
J. Appl. Probab.}, {\bf 21}, 768--801.

\item Gross, K.\,I.\ and Richards, D.\,St.\,P. (1998).
Algebraic methods toward higher-order probability inequalities.
In: {\it Stochastic Processes and Related Topics} (B. Rajput et
al., eds.), 189--211, Birkh\"{a}user, Boston.

\item Gross, K.\,I.\ and Richards, D.\,St.\,P. (2004).
Algebraic methods toward higher-order probability inequalities,
II. {\it Ann. Probab.}, {\bf 32}, 1509--1544.

\item Hoffmann-J$\o$rgensen, J. (2013). Slepian's inequality, modularity
and integral orderings. {\it High Dimensional Probability VI},
19--53, {Progress in Probability}, {66}, Springer, Basel.

\item Holland, P.\,W. and Wang, Y.\,J. (1987). Dependence function
for continuous bivariate densities. {\it Comm. Statist. -- Theory
and Methods}, {\bf 16}, 863--876.

\item Jones, M.\,C. (1996). The local dependence function. {\it Biometrika},
{\bf 83}, 899--904.

\item Jones, M.\,C. (1998). Constant local dependence.
{\it J. Multivariate Anal.}, {\bf 64}, 148--155.

\item Karlin, S. (1968). \textit{Total Positivity}, Vol.\,I. Stanford
University Press, CA.

\item Karlin, S. and Proschan, F. (1960). P\'olya type distributions of convolutions.
{\it Ann. Math. Statist.}, {\bf 31}, 721--736.

\item Kayid, M., Izadkhah, S. and Alshami, S. (2016). Laplace transform ordering of
time to failure in age replacement models. {\it J. Korean
Statist. Theory}, {\bf 45}, 101--113.

\item Kemperman, J.\,H.\,B. (1977). On the {\it FKG}--inequality for
measures on a partially ordered space. {\it Indagationes
Mathematicae}, {\bf 80}, 313--331.

\item Kotz, S., Balakrishnan, N. and Johnson, N.\,L. (2000). {\it
Continuous Multivariate Distributions}, Vol.\,1: Models and
Applications, 2nd ed. Wiley, New York.

\item {Kulkarni, H.\,V. (2006)}. Characterizations and modelling
of multivariate lack of memory property. {\it Metrika}, {\bf
64}, 167--180.

\item {Lai, C.-D. and Lin, G.\,D. (2014)}. Mean time to failure
of systems with dependent components. {\it Appl. Math. Comput.},
{\bf 246}, 103--111.

\item Lai, C.-D. and Xie, M. (2006). {\it Stochastic Ageing and Dependence for Reliability.} Springer, New York.

\item Lee, M.-L.T. (1985). Dependence by total positivity. {\it Ann. Probab.}, {\bf 13}, 572--582.

\item {Li, X. and Pellerey, F. (2011)}. Generalized
Marshall--Olkin distributions and related bivariate aging
properties. {\em J. Multivariate Anal.}, {\bf 102}, 1399--1409.

\item {Lin, G.\,D., Dou, X., Kuriki, S. and Huang, J.\,S.
(2014)}. Recent developments on the construction of bivariate
distributions with fixed marginals. {\it Journal of Statistical
Distributions and Applications}, {\bf 1}: 14.

\item {Lin, G.\,D., Lai, C.-D. and Govindaraju, K. (2016)}.
Correlation structure of the Marshall--Olkin bivariate exponential
distribution. {\it Statist. Methodology}, {\bf 29}, 1--9.

\item {Marshall, A.\,W. and Olkin, I. (1967)}.
 A multivariate exponential distribution. {\em J. Amer. Statist. Assoc.}, {\bf 62}, 30--44.

\item {Marshall, A.\,W. and Olkin, I. (2015)}. A bivariate
Gompertz--Makeham life distribution. {\em J. Multivariate
Anal.}, {\bf 139}, 219--226.

\item Marshall, A.\,W., Olkin, I.\ and Arnold, B.\,C. (2011).
\textit{Inequalities: Theory of Majorization and Its
Applications}, 2nd ed. Springer, New York.

\item {M\"uller, A. and Stoyan, D. (2002)}. {\it Comparison
Methods for Stochastic Models and Risks.} Wiley, New York.

\item Nadarajah, S. (2006). Exact distributions of $XY$ for some
bivariate exponential distributions. {\it Statistics}, {\bf 40},
307--324.

\item {Nelsen, R.\,B. (2006)}. {\it An Introduction to Copulas},
2nd ed. Springer, New York.

\item Newman, C. (1984). Asymptotic independence and limit theorems for positively and negatively
dependent random variables. In: {\it Inequalities in Statistics
and Probability}  (Y. L. Tong, ed.), IMS Lecture
Notes--Monograph Series, Vol. 5, 127--140.

\item {Rao, C.\,R. and Shanbhag, D.\,N. (1994)}. {\it
Choquet--Deny Type Functional Equations with Applications to
Stochastic Models.} Wiley, New York.

\item Shaked, M. (1977). A family of concepts of positive dependence for bivariate distributions.
{\it J. Amer. Statist. Assoc.}, {\bf 72}, 642--650.

\item {Shaked, M. and Shanthikumar, J.\,G. (1988)}. Multivariate conditional
hazard rates and the MIFRA and MIFR properties. {\it J. Appl.
Probab.}, {\bf 25}, 150--168.

\item {Shaked, M. and Shanthikumar, J.\,G. (2007)}. {\it
Stochastic Orders.} Springer, New Jersey.

\item Shorack, G.\,R. and Wellner, J.\,A. (1986).
{\it Empirical Processes with Applications to Statistics.}
Wiley, New York.

\item Slepian, D. (1962). The one-sided barrier problem for
Gaussian noise. {\it Bell System Technical Journal}, {\bf 41},
463-501.

\item Wang, Y.\,J. (1993). Construction of continuous bivariate density functions.
{\it Statist. Sinica}, {\bf 3}, 173--187.

\end{description}

\end{document}